\documentclass[11pt]{article} 
\usepackage{amsmath}
\usepackage{amsfonts}
\usepackage{amssymb}
\usepackage{graphicx}
\usepackage{color}
\usepackage[latin1]{inputenc}
\usepackage[T1]{fontenc}
\usepackage[french]{babel}

\makeatletter
\def\@sect#1#2#3#4#5#6[#7]#8{%
  \ifnum #2>\c@secnumdepth
    \let\@svsec\@empty
  \else
    \refstepcounter{#1}%
    \protected@edef\@svsec{\@seccntformat{#1}\relax}%
  \fi
  \@tempskipa #5\relax
  \ifdim \@tempskipa>\z@
    \begingroup
      #6{%
        \@hangfrom{\hskip #3\relax\@svsec}%
          \interlinepenalty \@M #8\@@par}%
    \endgroup
    \csname #1mark\endcsname{#7}%
    \addcontentsline{toc}{#1}{%
      \ifnum #2>\c@secnumdepth \else
        \protect\numberline{\csname the#1\endcsname.}%
      \fi
      #7}%
  \else
    \def\@svsechd{%
      #6{\hskip #3\relax
      \@svsec #8}%
      \csname #1mark\endcsname{#7}%
      \addcontentsline{toc}{#1}{%
        \ifnum #2>\c@secnumdepth \else
          \protect\numberline{\csname the#1\endcsname.}%
        \fi
        #7}}%
  \fi
  \@xsect{#5}}
\def\@seccntformat#1{\csname the#1\endcsname.\quad}
\makeatother

\newtheorem{theo}[equation]{Th\'eor\`eme}
\newtheorem{lem}[equation]{Lemme}
\newtheorem{proposition}[equation]{Proposition}
\newtheorem{definition}[equation]{D\'efinition}

\newtheorem{question}[equation]{Question}

\newenvironment{remarque}{
\refstepcounter{equation}\trivlist%
\item[\hskip \labelsep{\bfseries Remarque \theequation.\ }]}%
{\endtrivlist}%

\makeatletter
\renewcommand\theequation{\thesection.\arabic{equation}}
\@addtoreset{equation}{section}
\makeatother
\newcommand{\carrenoir}{\rule{0.5em}{0.5em}}
\makeatletter
\newenvironment{demo}[1][\@empty]{\textbf{D\'emonstration~%
\ifx\@empty#1:\else #1~:\fi~}}
{\hfill\carrenoir\nolinebreak\vspace{2mm}}
\makeatother

\newcommand{\oper}[2]{\newcommand{#1}{\mathop{\mathrm{#2}}\nolimits} }
\oper{\Vol}{Vol}
\newcommand{\R}{\mathbb R}
\newcommand{\de}{\mathrm{ d }}
\oper{\Ker}{Ker}
\oper{\Ima}{Im}

\DeclareSymbolFont{greek}{OML}{ptmcm}{m}{it}
\DeclareMathSymbol{\codiff}{\mathord}{greek}{"0E}
\DeclareMathSymbol{\prodint}{\mathord}{greek}{"13}
\newcommand{\N}{\mathbb N}

\title{Prescription de la multiplicité des valeurs propres du laplacien
de Hodge-de~Rham}
\author{Pierre Jammes}
\date{}
\begin{document}
\maketitle
{\small 
\textsc{Résumé.---}
Sur toute variété compacte de dimension supérieure ou égale à~6,
on prescrit le volume et le début du spectre du laplacien de Hodge-de~Rham
agissant sur les $p$-formes différentielles pour $1\leq p<\frac n2$.
En particulier, on prescrit la multiplicité des premières valeurs propres.

Mots-clefs : laplacien de Hodge-de~Rham, formes différentielles,
multiplicité de valeurs propres.

\medskip
\textsc{Abstract.---}
On any compact manifold of dimension greater than~6, we prescribe the volume
and any finite part of the spectrum of the Hodge Laplacian acting on $p$-form
for $1\leq p<\frac n2$. In particular, we prescribe the multiplicity 
of the first eigenvalues.


}

\section{Introduction}
On sait depuis les travaux de S.~Y.~Cheng \cite{ch76} que la multiplicité de la
$k$-ième valeur propre du laplacien sur une surface compacte est majorée
en fonction de $k$ et de la topologie. En dimension plus grande, 
Y.~Colin de Verdière a montré (\cite{cdv86}, \cite{cdv87}) que toute 
rigidité disparaît et qu'on peut arbitrairement prescrire le début du 
spectre, en particulier la multiplicité des valeurs propres peut être 
arbitrairement grande. 

Le résultat de Cheng s'étend aux opérateurs de Schrödinger sur les surfaces
et la majoration de la multiplicité a été améliorée (voir \cite{be80}, 
\cite{na88}, \cite{hohon99}), la meilleure estimation pour la multiplicité 
de la 2\ieme~valeur propre d'un opérateur de Schrödinger sur une surface
ayant été obtenue par B.~Sévennec (\cite{se94}, \cite{se02}). 
Ce problème a aussi été étudié pour des opérateurs avec
champ magnétique (\cite{cdvt93}, \cite{bcc98}, \cite{er02}), pour lesquels la 
multiplicité peut être arbitrairement grande. En dimension supérieure
ou égale à~3, J.~Lohkamp a amélioré les résultats de Colin de Verdière 
en montrant dans \cite{lo96} qu'on pouvait prescrire simultanément le 
début du spectre, le volume et certains invariants de courbure.

En comparaison, les connaissances sont beaucoup plus limités concernant
les opérateurs agissant sur les fibrés vectoriels naturels. P.~Guérini
a montré dans \cite{gu04} qu'on peut prescrire toute partie finie du 
spectre du laplacien de Hodge-de~Rham, qui agit sur les formes différentielles,
mais en imposant aux valeurs propres prescrites d'être simples. Un résultat 
semblable a été obtenu par M.~Dahl pour l'opérateur de Dirac (\cite{da05}).
Le seul résultat connu concernant la multiplicité des valeurs propres
non nulles de ces deux opérateurs est qu'on peut construire un nombre
arbitraire de valeurs propres doubles du laplacien de Hodge-de~Rham 
(voir \cite{ja09a}). On peut consulter \cite{ja09b} pour une présentation
générale de ces résultats.

Le but de cet article est d'étendre le théorème de Colin de Verdière aux 
formes différentielles en montrant que sur toute variété compacte de 
dimension $n\geq6$, on peut construire des valeurs propres du laplacien 
de Hodge-de~Rham de multiplicité arbitrairement grande, et plus précisément
que si on excepte les formes de degré $\frac n2$, on peut prescrire 
arbitrairement le début du spectre, avec multiplicité.

Si $(M^n,g)$ est une variété riemannienne
compacte orientable de dimension $n$, le laplacien $\Delta^p$ agissant sur
l'espace $\Omega^p(M)$ des $p$-formes différentielles est défini par
$\Delta=\de\codiff+\codiff\de$ où $\codiff$ désigne la codifférentielle.
Nous noterons son spectre
\begin{equation}
0=\lambda_{p,0}(M,g)<\lambda_{p,1}(M,g)\leq\lambda_{p,2}(M,g)\leq\ldots
\end{equation}
où les valeurs propres non nulles sont répétées s'il y a multiplicité.
La multiplicité de la valeur propre nulle, si elle existe, est un invariant
topologique: c'est le nombre de Betti $b_p(M)$. 
Par théorie de Hodge, le spectre $(\lambda_{p,i}(M,g))_{i\geq1}$ est la 
réunion de $(\mu_{p,i}(M,g))_i$ et $(\mu_{p-1,i}(M,g))_i$ où
\begin{equation}
0<\mu_{p,1}(M,g)\leq\mu_{p,2}(M,g)\leq\ldots
\end{equation} 
désigne les valeurs propres du laplacien restreint à l'espace des $p$-formes 
coexactes, et on a en outre $\mu_{p,i}(M,g)=\mu_{n-p-1,i}(M,g)$ pour tout 
$p$ et $i$ si $M$ n'a pas de bord. Le spectre complet du laplacien se 
déduit alors des $\mu_{p,i}(M,g)$ pour $p\leq\frac{n-1}2$, c'est donc
à la multiplicité de ces valeurs propres qu'on va s'intéresser.

\begin{theo}\label{intro:th1}
Soit $M^n$ une variété compacte connexe orientable
sans bord de dimension $n\geq6$ et $N\in\N^*$.
Si on se donne un réel $V>0$, une suite 
$0<a_{1,1}<a_{1,2}\leq a_{1,3}\leq\ldots\leq a_{1,N}$
et des suites $0<a_{p,1}\leq a_{p,2}\leq\ldots\leq a_{p,N}$
pour $2\leq p\leq[\frac{n-3}2]$, alors il existe une métrique $g$ sur $M$ 
telle que
\begin{itemize}
\item $\mu_{p,k}(M,g)=a_{p,k}$ pour $1\leq k\leq N$ et 
$1\leq p\leq[\frac{n-3}2]$;
\item $\mu_{[\frac{n-1}2],1}(M,g)>\sup_{i,N}a_{p,i}$;
\item $\Vol(M,g)=V$.
\end{itemize}
\end{theo}
\begin{remarque}
La condition $\mu_{[\frac{n-1}2],1}(M,g)>\sup_{i,N}a_{p,i}$ permet
de prescrire le début du spectre en degré $[\frac{n-1}2]$, les formes 
propres correspondantes étant exactes.
\end{remarque}
 Comme dans \cite{cdv86} et \cite{cdv87}, le principe de la démonstration
consiste à faire converger le spectre de la variété vers celui d'un espace 
modèle (en l'occurence un domaine de la variété) qui a le spectre souhaité 
et de conclure grâce aux propriétés de stabilité du spectre du modèle.
 Pour les formes de degrés proches de $\frac n2$, la démonstration
échoue en raison de l'invariance (ou la presque invariance) conforme
de la norme $L^2$. On verra cependant que les rigidités qui apparaissent
pour ces degrés faciliteront la construction des espaces modèles.
Une autre difficulté apparaît pour les formes de degré~1, on ne pourra
pas faire appel aux même espaces modèles que pour les autres degrés. 
On utilisera une construction particulière qui ne permet pas de prescrire
la multiplicité de la première valeur propre et qui ne fonctionne qu'en
dimension $n\geq6$.

 Le problème suivant reste donc ouvert :
\begin{question}
La multiplicité des valeurs propres $\mu_{1,1}(M,g)$ et 
$\mu_{[\frac{n-1}2],k}(M,g)$ peut-elle être arbitrairement grande ?
\end{question}
Dans \cite{ja09a}, on montre comment construire des exemples de variétés 
de dimension $n\geq4$ admettant des valeurs propres de multiplicité 
arbitrairement grande, y compris en degré $[\frac{n-1}2]$. Leur topologie 
est très particulière (variétés produits), mais ces exemples montrent
qu'on a pas en général de borne sur la multiplicité comme en dimension~2.

 En ce qui concerne $\mu_{1,1}(M,g)$, on peut aussi apporter cet élément 
de réponse en utilisant les même techniques que pour le 
théorème~\ref{intro:th1}:
\begin{theo}\label{intro:th2}
Si $M^n$ une variété compacte connexe orientable sans bord de dimension
$n\geq5$, alors il existe sur $M$ une métrique $g$ telle que la multiplicité
de $\mu_{1,1}(M,g)$ soit égale à~3.
\end{theo}
La méthode utilisée ne permet cependant pas de prescrire les autres valeurs
propres.

 Le problème le plus intéressant semble être de comprendre ce qui se passe 
en dimension~3. En effet, les exemples produits donné dans \cite{ja09a}
sont de dimension au moins~4. L'énoncé de S.~Y.~Cheng pourrait s'étendre
aux 1-formes coexactes en dimension~3:
\begin{question}
Sur une variété $M$ de dimension~3, existe-t-il une borne sur la multiplicité 
de $\mu_{1,k}(M,g)$ dépendant uniquement de $k$ et de la topologie de~M ?
\end{question}
Il faut noter que pour établir un tel résultat, on peut difficilement
espérer adapter de la démonstration de Cheng dont les arguments sont
spécifiques aux fonctions (domaines nodaux) et à la topologie en dimension~2.

 Dans la section~\ref{conv} nous montrerons, après quelques rappels 
techniques, qu'on peut faire tendre le spectre d'une variété compacte pour
le laplacien de Hodge-de~Rham vers le spectre d'un de ses domaines. La 
section~\ref{presc} sera consacrée à la démonstration des 
théorèmes~\ref{intro:th1} et \ref{intro:th2}.

\section{Convergence du spectre d'une variété vers celui d'un domaine}
\label{conv}
\subsection{Conditions de bord et cohomologie pour une variété à bord}
\label{conv:bord}
Dans ce paragraphe et le suivant, nous allons rappeler certains aspects
techniques de la théorie spectrale du laplacien de Hodge-de~Rham auquels
nous feront appel pour montrer qu'on peut faire tendre le spectre d'une
variété compacte vers celui d'un de ses domaines (théorème~\ref{conv:th}).

Si $U$ est un domaine à bord $C^1$ d'une variété compacte $M$, 
on note $j:\partial U\to \overline U$ 
l'injection canonique et $N$ un champ de vecteur normal au bord.
Il existe plusieurs conditions de bord admissibles pour le laplacien de
Hodge-de~Rham sur $U$ (c'est-à-dire telles que le laplacien soit elliptique), 
les deux principales sont les conditions absolues et
relatives, que nous noterons respectivement (A) et (R) et qui sont
définies par
\begin{equation}
(A)\ \left\{\begin{array}{l}j^*(\prodint_N\omega)=0\\
j^*(\prodint_N\de\omega)=0\end{array}\right.\textrm{ ou }
\left\{\begin{array}{l}j^*(*\omega)=0\\
j^*(*\de\omega)=0\end{array}\right.
\end{equation}
et
\begin{equation}
(R)\ \left\{\begin{array}{l}j^*(\omega)=0\\
j^*(\codiff\omega)=0\end{array}\right.
\end{equation}
Pour la condition (A), $\Ker\Delta$ est isomorphe à la cohomologie $H^p(U)$
et pour (R), il est isomorphe à la cohomologie à support compact $H^p_0(U)$
(voir par exemple \cite{ta96}, ch.~5). Rappelons aussi que les cohomologies 
de $U$ et $\overline U$ sont isomorphes (\cite{ta96} ch.~5 p.~375). Il est 
immédiat que sous la 
condition (R) on a $j^*(\de\omega)=0$. Et comme la dualité de Hodge
permute ces deux conditions de bord, (A) implique que $j^*(*\codiff\omega)=0$.

 Nous aurons besoin d'une autre condition de bord, définie par 
\begin{equation}
(D)\ \left\{\begin{array}{l}j^*(\omega)=0\\
j^*(*\omega)=0\end{array}\right.
\end{equation}
Pour la condition (D), le noyau du laplacien est trivial (voir~\cite{an89}).

 Rappelons qu'en restriction aux fonctions, la condition (A) est 
équivalente à la condition de Neumann, et les conditions (R) et (D) 
à la condition de Dirichlet.

 La décomposition de Hodge $L^2(\Lambda^p\overline U)=\Ima\de\oplus\Ker\Delta
\oplus\Ima\codiff$ dépend de la condition de bord choisie. Pour une forme
$\omega$, elle s'écrit $\omega=\de\codiff\theta+\alpha+\codiff\de\theta$
où $\alpha$ et $\theta$ vérifie la condition de bord considérée, $\alpha$
étant dans le noyau $\Ker\Delta$ correspondant (cf. proposition~9.8 de
\cite{ta96}).

 Dans la suite du texte interviendra essentiellement la 
condition~(A) (et dans une moindre mesure la condition~(D)), et nous
ferons en particulier appel à la propriété suivante :
\begin{proposition}\label{conv:lem:prim}
Si $\omega$ est une $p$-forme exacte sur $\overline U$, alors il existe une
forme $\theta\in\Omega^p(\overline U)$ vérifiant la condition (A) et telle que 
$\omega=\de\codiff\theta$. En particulier, la forme $\varphi=\codiff\theta$ 
vérifie $j^*(\prodint_N\varphi)=0$, est orthogonale aux formes fermées 
et minimise la norme $L^2$ parmi les primitives de $\omega$.
\end{proposition}
 La démonstration de ces faits (qui ne sont souvent utilisés
qu'implicitement dans la littérature) est esquissée dans \cite{ch79} 
(section~3, p.~272). Nous allons la rappeler :

\begin{demo}
Si $\alpha\in\Omega^p(\overline U)$ et $\beta\in\Omega^{p+1}(\overline U)$, 
une intégration par partie donne $(\de\alpha,\beta)=(\alpha,\codiff\beta)+
\int_{\partial U}j^*(\alpha\wedge*\beta)$. Pour que le terme de bord 
s'annule, il n'est pas nécessaire que les deux formes $\alpha$ et $\beta$ 
vérifient l'une des conditions de bord, il suffit par exemple que 
$j^*(*\beta)=0$, c'est-à-dire que $\beta$ soit tangentielle.
 
 Si $\omega$ est une forme exacte, on a alors $(\omega,\beta)=0$ pour
toute forme $\beta$ tangentielle cofermée. Dans la décomposition 
de Hodge pour la condition~(A), la composante cofermée de $\omega$ est
donc nulle même si $\omega$ ne vérifie aucune condition de bord. 
On a donc $\omega=\de\codiff\theta$ où $\theta$ vérifie la condition~(A)
(autrement dit, l'adhérence $L^2$ des formes exactes vérifiant (A) contient
toutes les formes exactes, sans condition de bord).
 
 Si on pose $\varphi=\codiff\theta$, on a alors $j^*(\prodint_N\varphi)=0$, 
et $(\alpha,\varphi)=(\de\alpha,\theta)=0$ pour toute forme $\alpha$ fermée.
Comme deux primitives de $\omega$ diffèrent par une forme fermée,
$\varphi$ minimise la norme $L^2$ parmi celles-ci.
\end{demo}

 Nous aurons besoin
de définir, outre les cohomologies $H^p(U)$ et $H^p_0(U)$ évoquées
plus haut, un espace de cohomologie traduisant l'interaction entre
la cohomologie de $U$ et celle de $M$. Cet espace est construit comme 
le quotient des formes fermées de $U$ par la restriction des formes
fermées de $M$ : 
\begin{equation}
H^p(U/M)=\{\omega\in\Omega^p(\overline U),\ \de\omega=0\}/
\{\omega_{|\overline U},\ \omega\in\Omega^p(M)\textrm{ et }\de\omega=0\}
\end{equation}
 Comme une forme exacte de $\overline U$ est toujours la restriction d'une 
forme exacte de $M$, $H^p(U/M)$ est isomorphe au quotient de $H^p(U)$ par 
l'image de l'application naturelle $H^p(M)\to H^p(U)$ définie par restriction
des formes fermées et exactes. En particulier, $H^p(U/M)$ est de
dimension finie.

\subsection{Caractérisation du spectre du laplacien de Hodge}
Pour démontrer le théorème de convergence au paragraphe suivant, nous
utiliserons une caractérisation variationnelle du spectre dont le
principe est dû à J.~Cheeger et J.~Dodziuk :
\begin{proposition}[\cite{do82}, \cite{mc93}]\label{conv:prop1}
 Sur une variété compacte sans bord ou avec condition de bord (A), on a
$$\mu_{p,i}=\inf_{V_i}\sup_{\omega\in V_i\backslash\{0\}}\left\{
\frac{\|\omega\|^2}{\|\varphi\|^2},\ \de\varphi=\omega\right\},$$
où $V_i$ parcourt l'ensemble des sous-espaces de dimension $i$ dans
l'espace des $p+1$-formes exactes lisses.
\end{proposition}
Il faut préciser que dans le cas à bord, cette formule fournit le spectre 
pour la condition (A) même si on ne suppose pas que des formes $\omega$ et 
$\varphi$ vérifient cette condition. Cela tient essentiellement au fait 
démontré plus haut que l'espace $L^2$ des formes exactes pour la 
décomposition de Hodge associée à (A) contient toutes les formes
exactes sans condition de bord.

 La proposition~\ref{conv:prop1} permet d'estimer les valeurs propres
du laplacien, mais pour les applications à la multiplicité, et en particulier 
pour utiliser les techniques de Colin de Verdière (voir 
lemmes~\ref{conv:lem1} et~\ref{conv:lem2}), nous aurons
besoin de contrôler à la fois les valeurs propres et les espaces propres,
nous allons donc la reformuler. Si $\omega$ est 
une forme exacte, alors $q(\omega)=\inf_{\de\varphi=\omega}\|\varphi\|^2$ 
est une forme quadratique, c'est la norme au carré de la primitive coexacte 
de $\omega$. Son spectre est l'inverse de celui du laplacien (on peut
se convaincre qu'on a aussi $q(\omega)=(\Delta^{-1}\omega,\omega)$). 
On retrouve la formule de la proposition~\ref{conv:prop1} en intervertissant
le rôle de la forme quadratique et de la norme de Hilbert :
\begin{proposition}\label{conv:prop2}
Le spectre et les espaces propres du laplacien en restriction aux formes
exactes sont ceux de la forme quadratique $Q(\omega)=\|\omega\|^2_{L^2}$
relativement à la norme $|\omega|=\displaystyle\inf_{\de\varphi=\omega}
\|\varphi\|_{L^2}$.
\end{proposition}
\begin{remarque}
L'espace $\overline{L^2(\Lambda^{p+1}M)\cap\Ima\de}$ (l'adhérence étant
prise au sens de la norme $L^2$) n'est pas un espace 
de Hilbert car il n'est pas complet pour la norme $|\cdot|$ (c'est 
seulement le domaine de la forme quadratique $Q$), mais on peut identifier 
son complété à $L^2(\Lambda^pM)/\overline{\Ker\de}$, chaque forme exacte de 
$L^2(\Lambda^{p+1}M)$ étant identifiée à l'ensemble de ses primitives.
\end{remarque}
\begin{remarque}
On peut déduire aisément de la proposition~\ref{conv:prop1} que le spectre 
du laplacien est continu pour la topologie $C^0$ (cf.~\cite{do82}). 
On sait aussi, d'après \cite{bd97} qu'on a continuité des espaces propres 
pour les topologies $C^0$ et $L^2$. La formulation de la 
proposition~\ref{conv:prop2} permet de retrouver ce fait de manière plus 
directe à l'aide des techniques de \cite{cdv86} (voir le 
lemme~\ref{conv:lem2} au paragraphe suivant). 
\end{remarque}

\subsection{Convergence spectrale}
Nous allons maintenant montrer comment on peut faire tendre
le spectre du laplacien de Hodge-de~Rham d'une variété compacte vers
celui d'un domaine, généralisant ainsi au formes différentielles le
théorème que Y.~Colin de Verdière a montré pour les fonctions (théorème~III.1 
de \cite{cdv86}). Pour appliquer ce résultat, nous aurons besoin d'une 
certaine uniformité de la convergence, nous reprendrons pour cela les
notations de \cite{cdv86} :

 Soit $E_0$ et $E_1$ sont deux sous-espaces vectoriels de même dimension~$N$
d'un espace de Hilbert, munis respectivement des formes quadratiques $q_0$
et $q_1$. Si $E_0$ et $E_1$ sont suffisamment proches, il existe
une isométrie naturelle $\psi$ entre les deux (voir la section~I de
\cite{cdv86} pour les détails de la construction), on définit alors
l'écart entre $q_0$ et $q_1$ par $\|q_1\circ\psi-q_0\|$. Pour deux
formes quadratiques $Q_0$ et $Q_1$ sur l'espace de Hilbert, on appellera
\emph{$N$-écart spectral entre $Q_0$ et $Q_1$} l'écart entre les
deux formes quadratiques restreintes à la somme des espaces propres associés
aux $N$ premières valeurs propres. Si cet écart est petit, alors les $N$ 
premières valeurs propres de $Q_0$ et leurs espaces propres sont proches 
de ceux de $Q_1$. 

 On veut montrer que la convergence spectrale est uniforme pour une certaine
famille de spectres limites. Comme dans \cite{cdv86} on dira donc qu'une forme
quadratique vérifie l'hypothèse ($*$) si ses valeurs propres vérifient
$$\lambda_1\leq\ldots\leq\lambda_N<\lambda_N+\eta\leq\lambda_{N+1}\leq M$$
pour un entier $N$ et des réels $\eta,M>0$ fixés une fois pour toute.

 Dans la suite du texte, sauf mention explicitement contraire, le spectre 
considéré sur les domaines sera toujours relatif aux conditions de bord~(A).
\begin{theo}\label{conv:th} 
Soit $(M^n,g)$ une variété riemanienne compacte sans bord de dimension $n$
et $U$ un domaine de $M$ à bord $C^1$. Il existe une suite de
métriques $(g_i)$ sur $M$ conformes à $g$ telle que 
\begin{itemize}
\item $\Vol(M,g_i)\to\Vol(U,g)$ quand $i\to\infty$;
\item $\mu_{p,k}(M,g_i)\to0$ pour $p\leq\frac{[n-3]}2$ et $k\leq d_p$ quand
$i\to\infty$;
\item $\mu_{p,k+d_p}(M,g_i)\to\mu_{p,k}(U,g)$ pour $p\leq\frac{[n-3]}2$ et 
$k\geq1$ quand $i\to\infty$;
\end{itemize}
où $d_p$ est la dimension de $H^p(U/M)$.

En outre, si les $\mu_{p,k}(U,g)$ vérifient l'hypothèse ($*$) pour
$1\leq p\leq\frac{[n-3]}2$, alors pour tout $\varepsilon>0$ il existe 
$i$ tel que le $N$-écart spectral entre les laplaciens sur $U$ et $M$
pour la métrique $g_i$ soit inférieur à $\varepsilon$.
\end{theo}
\begin{remarque} 
Dans \cite{co04}, B.~Colbois avait posé la question de savoir si on 
peut faire tendre les valeurs propres du laplacien de Hodge-de~Rham 
vers~0 en fixant le volume et la classe conforme. Pour $p\leq\frac{[n-3]}2$,
une réponse positive a été donnée dans \cite{ja08} par une constrution
similaire à celle du théorème~\ref{conv:th} avec $U\simeq
S^p\times B^{n-p}$. Le théorème~\ref{conv:th} permet une compréhension
plus générale de ce phénomène, tout en simplifiant les démonstrations
de \cite{ja08}.
\end{remarque}
On sait que l'énoncé du théorème~\ref{conv:th} est faux pour 
$p=[\frac{n-1}2]$, car pour ce degré, le spectre du laplacien est
uniformément minoré dans une classe conforme:
\begin{theo}[\cite{ja07}]\label{conv:th2}
Soit $M^n$ une variété compacte de dimension $n\geq3$, Pour toute
classe conforme $C$ sur $M$, il existe une constante $K(C)>0$ telle que pour
toute métrique $g\in C$, on a
$$\mu_{\left[\frac{n-1}2\right],1}(M,g)\Vol(M,g)^{\frac2n}\geq 
K.$$
\end{theo}
 On peut donc construire des contre-exemples en choissant $M$ et $U$ tels
que $H^{[\frac{n-1}2]}(M/U)$ soit non trivial, ou en utilisant un domaine
dont la première valeur propre est plus petite que que 
$K\cdot\Vol(U,g)^{-\frac2n}$ (c'est possible d'après~\cite{gu04}).
\begin{remarque}\label{conv:rq2}
La constante $K$ du théorème~\ref{conv:th2} varie continûment pour des 
déformations $C^0$ de la métrique (voir~\cite{ja07}). Cette propriété nous
sera utile pour démontrer le théorème~\ref{intro:th1}.
\end{remarque}
\begin{remarque}
Pour l'opérateur de Dirac, il existe une rigidité conforme du même type
que celle du théorème~\ref{conv:th2} (voir~\cite{am03}). On ne peut donc pas 
espérer obtenir un résultat semblable au  théorème~\ref{conv:th} pour 
les spineurs.
\end{remarque}

Comme dans \cite{cdv86}, on fera appel aux deux lemmes qui suivent. Les
constantes $N$, $M$ et $\eta$ qui interviennent dans les énoncés
font référence à l'hypothèse ($*$) définie plus haut.
\begin{lem}[\cite{cdv86}, th.~I.7]\label{conv:lem1}
Soit $Q$ une forme quadratique positive sur un espace de Hilbert $\mathcal H$ 
dont le domaine admet la décomposition $Q$-orthogonale 
$\mathrm{dom}(Q)=\mathcal H_0\oplus\mathcal H_\infty$. Pour tout 
$\varepsilon>0$, il existe une constante $C(\eta,M,N,\varepsilon)>0$ (grande) 
telle que si $Q_0=Q_{|\mathcal H_0}$ vérifie l'hypothèse ($*$) et que 
$\forall x\in\mathcal H_\infty,\ Q(x)\geq C|x|^2$,
alors $Q$ et $Q_0$ ont un $N$-écart spectral inférieur à $\varepsilon$.
\end{lem}

\begin{lem}[\cite{cdv86}, th.~I.8]\label{conv:lem2}
Soit $(\mathcal H,|\cdot|)$ un espace de Hilbert muni d'une forme quadratique
positive $Q$. On se donne en outre une suite de métriques $|\cdot|_n$
sur $\mathcal H$ et une suite de formes quadratiques $Q_n$ de même
domaine que $Q$ telles que:
\begin{itemize}
\item[(i)] il existe $C_1,C_2>0$ tels que $\forall x\in\mathcal H,\
C_1|x|\leq |x|_n\leq C_2|x|$;
\item[(ii)] pour tout $x\in \mathrm{dom}(Q)$, $|x|_n\to|x|$;
\item[(iii)] pour tout $x\in \mathrm{dom}(Q)$, $Q(x)\leq Q_n(x)$;
\item[(iv)] pour tout $x\in \mathrm{dom}(Q)$, $Q_n(x)\to Q(x)$.
\end{itemize}
Si $Q$ vérifie l'hypothèse ($*$), alors à partir d'un certain rang (dépendant
de $\eta$, $M$ et $N$), $Q$ et $Q_n$ ont un $N$-écart spectral inférieur à
$\varepsilon$.
\end{lem}
\begin{remarque}\label{conv:rem1}
Dans le lemme~\ref{conv:lem2}, on peut affaiblir l'hypothèse
$C_1|x|\leq |x|_n\leq C_2|x|$ en $C_1|x|\leq |x|_n\leq C_2|x|+
\varepsilon_nQ(x)^\frac12$ avec $\varepsilon_n\to0$,
la démonstration restant exactement la même. En particulier, il n'est
pas nécessaire que l'espace de Hilbert $(\mathcal H,|\cdot|)$ soit
complet pour $|\cdot|_n$.
\end{remarque}
\begin{remarque}
Pour déduire la convergence du spectre et des espaces propres de la 
convergence des formes quadratiques, on doit en principe se ramener à une
norme de Hilbert fixe. Ça ne sera pas nécessaire dans la suite car les
étapes de la démonstration où la norme varie seront traitées à l'aide
du lemme~\ref{conv:lem2}.
\end{remarque}
\begin{demo}[du théorème~\ref{conv:th}]
 En vertu de la proposition~\ref{conv:prop2}, il suffit de montrer
la convergence des valeurs propres et des espaces propres de la
forme quadratique $Q(\omega)=\|\omega\|^2$ relativement à la norme
$|\omega|=\inf_{\de\varphi=\omega}\|\varphi\|$ pour les formes 
différentielles exactes de degré 1 à $[\frac{n-1}2]$ (les $\mu_{p,i}$
sont les valeurs propres de $Q$ sur les formes exactes de degré $p+1$).

Comme dans \cite{cdv86}, on va passer par l'intermédiaire d'une famille
de  métriques singulières $(g_\varepsilon)_{\varepsilon\in]0,1]}$ définies
par $g_\varepsilon=g$ sur $U$ et $g_\varepsilon=\varepsilon^2g$
sur $M\backslash U$. Notons que la forme quadratique $Q_\varepsilon$ et
la métrique $|\cdot|_\varepsilon$ associées à $g_\varepsilon$ sont bien
définies. Plus précisément, pour une forme exacte $\omega\in\Omega^{p+1}(M)$,
on a
\begin{equation}
Q_\varepsilon(\omega)=\int_U|\omega|^2\de v_g+\varepsilon^{n-2p-2}
\int_{M\backslash U}|\omega|^2\de v_g
\end{equation}
et
\begin{equation}
|\omega|_\varepsilon^2=\inf_{\de\varphi=\omega}\left(\int_U|\varphi|^2\de v_g
+\varepsilon^{n-2p}\int_{M\backslash U}|\varphi|^2\de v_g\right).
\end{equation}
Pour $\varepsilon$ fixé, les normes induites par $g$ et $g_\varepsilon$ sont
équivalentes. Il en va donc de même pour les normes $|\cdot|$ et
$|\cdot|_\varepsilon$, ainsi que pour les normes d'opérateurs associées à
$Q$ et $Q_\varepsilon$. Les résultats usuels de théorie spectrale
(spectre formé d'une suite de valeurs propres tendant vers l'infini)
s'appliquent donc à la forme quadratique $Q_\varepsilon$ .

 La démonstration se déroule en trois étapes. D'abord, on montre 
que pour $\varepsilon$ donné, on peut approcher la métrique $g_\varepsilon$ 
par une métrique lisse avec convergence du spectre. Comme pour tout 
$\varepsilon$ on peut trouver une métrique lisse dont l'écart spectral avec 
$g_\varepsilon$ soit arbitrairement petit, il suffit de montrer que 
l'écart spectral entre le spectre de $g_\varepsilon$ et celui de $U$ 
devient lui ausssi arbitrairement petit quand $\varepsilon\to0$, ce qui
fait l'objet de la suite de la démonstration. 
La deuxième étape consiste à décomposer l'espace
$\mathcal H$ des formes exactes en une somme $\mathcal H_0\oplus
\mathcal H_\infty$ à laquelle on applique le lemme~\ref{conv:lem1}. Enfin,
on montre la convergence du spectre de $Q_\varepsilon$ restreint à 
$\mathcal H_0$ vers le spectre du domaine à l'aide du lemme~\ref{conv:lem2}.

\emph{Étape~1.}
En utilisant le lemme~\ref{conv:lem2}, on va montrer qu'on peut approcher 
$g_\varepsilon$ par des metriques lisses avec convergence du volume, 
du spectre et des espaces propres.

 Pour un $\varepsilon$ fixé, on peut approcher la fonction
$\chi_U+\varepsilon\chi_{M\backslash U}$ par une suite de fonctions
décroissantes $(f_j)$. La suite de métriques $g_j=f_j^2\cdot g$ tend alors vers
$g_\varepsilon$, et on note $Q_j$ et $|\cdot|_j$ la forme quadratique et
la norme hilbertienne associées. On peut vérifier que la suite de forme
quadratique $Q_j$ converge simplement vers $Q_\varepsilon$ (hypothèse~(iv)
du lemme~\ref{conv:lem2}). Comme il existe une constante $C$ telle que
$g_\varepsilon\leq g_j\leq C\cdot g_\varepsilon$, on a
$|\omega|_\varepsilon^2\leq |\omega|_j^2\leq C\cdot |\omega|_\varepsilon^2$
et $Q_\varepsilon(\omega)\leq Q_j(\omega)$ pour toute forme exacte $\omega$
(hypothèses~(i) et~(iii) du lemme). Pour appliquer le lemme~\ref{conv:lem2},
il reste à vérifier que l'hypothèse~(ii) est satisfaite, à savoir que
$|\cdot|_j$ converge simplement vers $|\cdot|_\varepsilon$. Étant donné
$\eta>0$, il existe une forme $\varphi_0$ telle que $\de\varphi_0=\omega$
et $\|\varphi_0\|_{g_\varepsilon}^2\leq
|\omega|_\varepsilon^2+\eta$ (par définition de $|\cdot|_\varepsilon$). Comme
$g_j$ converge vers $ g_\varepsilon$, pour $j$ assez grand on a aussi
$\|\varphi_0\|_{g_j}^2\leq\|\varphi_0\|_{g_\varepsilon}^2+\eta$, et donc
$|\omega|_\varepsilon^2\leq|\omega|_j^2\leq \|\varphi_0\|_{g_j}^2\leq
|\omega|_\varepsilon^2+2\eta$. On a donc bien $|\omega|_j
\to|\omega|_\varepsilon$ pour tout $\omega$.

 Selon le lemme~\ref{conv:lem2}, à $\varepsilon$ fixé, on peut donc trouver
un $j_\varepsilon$ tel que l'écart spectral entre $Q_{j_\varepsilon}$ et
$Q_\varepsilon$ soit arbitrairement petit. 

\emph{Étape~2.}
On va décomposer l'espace des $(p+1)$-formes exactes en une somme
$\mathcal H_0\oplus\mathcal H_\infty$ à laquelle on applique le
lemme~\ref{conv:lem1} pour la norme $|\cdot|_\varepsilon$ et la forme
quadratique~$Q_\varepsilon$.

On commence par définir le sous-espace $\mathcal H_\infty$ de 
$\overline{\Ima\de^p}\subset L^2(\Lambda^{p+1}M)$ comme l'adhérence des 
différentielles des formes lisses qui s'annulent sur $U$ (les adhérences
sont au sens de la norme $L^2$). Une telle forme 
va nécessairement vérifier la condition de bord de Dirichlet sur 
$M\backslash U$:
\begin{equation}
\mathcal H_\infty=\overline{\{\de\varphi,\ \varphi\in\Omega^p(M),\ 
\varphi_{|U}=0,\ \varphi_{|M\backslash U}\textrm{ vérifie }(D)\}}.
\end{equation}

L'espace $\mathcal H_0$ sera défini comme la somme de deux espaces
$\mathcal H_1$ et $\mathcal H_2$ construits séparément.

 Soit $\omega$ une $(p+1)$-forme exacte sur $\overline U$ et 
$\varphi\in\Omega^p(\overline U)$
la primitive coexacte de  $\omega$ fournie par le lemme~\ref{conv:lem:prim}. 
Si $\tilde\varphi$ est un prolongement lisse
de $\varphi$ sur $M$, alors $\de\tilde\varphi$ est définie à un élement
de $\mathcal H_\infty$ près. On peut définir alors $\tilde\omega$ comme
le $\de\tilde\varphi$ de norme minimale pour $g_\varepsilon$. Cet infimum
est bien atteint dans $L^2$ et on peut le construire par
projection sur l'orthogonal $L^2$ de $\mathcal H_\infty$. On pose alors
\begin{equation}
\mathcal H_1=\overline{\{\tilde\omega,\ \omega\in\Omega^{p+1}(\overline U)
\textrm{ exacte}\}}.
\end{equation}

 L'espace $\mathcal H_2$ est défini à partir de l'espace de cohomologie
$H^p(U/M)$. Comme $H^p(U/M)$ est isomorphe au quotient de $H^p(U)$ par le
sous-espace induit par $H^p(M)$, il est aussi isomorphe à l'orthogonal, 
dans l'espace des formes harmoniques de $U$ (avec condition de bord absolue),
des représentants harmoniques des classes de cohomologie induites par
$H^p(M)$. On peut ainsi définir un représentant harmonique $h$
de chaque classe $[c]\in H^p(U/M)$, qui est la forme $h\in[c]$ orthogonale 
aux restrictions des formes fermées de $M$.

On peut alors construire $\mathcal H_2$ sur le modèle de $\mathcal H_1$.
Chaque forme harmonique $h$ représentant une classe de $H^p(U/M)$,
peut être étendue en une forme $\tilde h$ sur $M$, la forme $\de\tilde h$
étant alors définie à un élement de $\mathcal H_\infty$ près. On notant
$\tilde\omega_h$ la forme $\de\tilde h$ qui est $L^2$-orthogonale à 
$\mathcal H_\infty$, on pose
\begin{equation}
\mathcal H_2=\{\tilde\omega_h,\ [h]\in H^p(U/M)\}.
\end{equation}

Par construction, $\mathcal H_0=\mathcal H_1\oplus\mathcal H_2$ et
$\mathcal H_\infty$ sont orthogonaux pour la norme $L^2$, donc $Q$-orthogonaux. 
Avant d'appliquer le lemme~\ref{conv:lem1} on doit encore vérifier que 
$\mathcal H_0\oplus\mathcal H_\infty$ contient bien toutes les 
$(p+1)$-formes exactes. Soit $\omega\in\Omega^{p+1}$ une forme exacte. Par
définition de $\mathcal H_1$, on peut écrire $\omega=\omega_1+\omega'$, avec
$\omega_1\in\mathcal H_1$ telle que $\omega_{|U}={\omega_1}_{|U}$ et 
$\omega'_{|U}=0$. Comme $\omega$ et $\omega_1$
sont exacte, $\omega'$ l'est aussi. Si on pose $\omega'=\de\varphi$, la
forme $\varphi$ est définie à une forme fermée près et vérifie
$\de\varphi=0$ sur $U$, la classe $[\varphi]\in H^p(U/M)$ est donc
bien définie. Par définition de $\mathcal H_2$ on a alors $\omega'=
\tilde\omega_h+\omega_0$, où $\tilde\omega_h$ est l'élément de $\mathcal H_2$
associé au représentant harmonique $h$ de $[\varphi]\in H^p(U/M)$, et
$\omega_0\in\mathcal H_\infty$. On a donc bien $\omega\in\mathcal H_1
\oplus\mathcal H_2\oplus\mathcal H_\infty$, et donc 
$\mathrm{dom}(Q_\varepsilon)=\mathcal H_0\oplus\mathcal H_\infty$.

Si on note $\lambda^{(D)}$ la première valeur propre de $M\backslash U$ pour
la métrique $g$ et la condition de bord (D), on sait que pour toute forme 
$\omega\in \mathcal H_\infty$ il existe, par définition, une forme 
$\varphi\in\Omega^p(M)$ à support dans $M\backslash U$ telle que 
$\de\varphi=\omega$ et $\|\omega\|^2/\|\varphi\|^2\geq\lambda^{(D)}$ pour
la métrique $g$. Pour la métrique $g_\varepsilon$, on a donc
$\|\omega\|^2/\|\varphi\|^2\geq\varepsilon^{-2}\lambda^{(D)}$, et \emph{a
fortiori} $Q_\varepsilon(\omega)/|\omega|_\varepsilon\geq
\varepsilon^{-2}\lambda^{(D)}$. Si $\varepsilon$ est suffisamment petit, on 
peut appliquer le lemme~\ref{conv:lem1}.

\emph{Étape~3.}
 On va achever la démonstration en appliquant le lemme~\ref{conv:lem2} et 
la remarque~\ref{conv:rem1} à l'espace $\mathcal H_0$ et aux familles 
de métriques et de formes quadratiques $|\cdot|_\varepsilon$ et 
$Q_\varepsilon$. On définit la forme quadratique $Q$ sur $\mathcal H_0$
par $Q(\omega)=\int_U|\omega|^2\de v_g$. On doit aussi définir une
norme $|\cdot|$ sur $\mathcal H_0$ : pour $\omega\in\mathcal H_1$,
on note $\varphi_\omega$ la primitive coexacte de $\omega_{|\overline U}$ 
donnée par la
proposition~\ref{conv:lem:prim}, et pour $\omega\in\mathcal H_2$, on note
$\varphi_\omega$ le représentant harmonique de la classe de cohomologie
définie par $\omega$. On étend linéairement l'application 
$\omega\to\varphi_\omega$ et on pose $|\omega|=\|\varphi_\omega\|$, la
norme $\|\cdot\|$ étant ici la norme $L^2$ sur les $p$-formes de~$U$.
Les espaces $\mathcal H_1$ et $\mathcal H_2$ sont orthogonaux pour
$\|\cdot\|$, le noyau de la forme quadratique $Q$ (relativement à $\|\cdot\|$) 
est $\mathcal H_2$, et le spectre de $Q$ sur $\mathcal H_1$ est le spectre
du domaine $U$ pour les $(p+1)$-formes exactes.

 On doit maintenant vérifier que les quatre hypothèses du 
lemme~\ref{conv:lem2} sont satisfaites.

Les hypothèses~(iii) et~(iv) sont les plus simples à vérifier~: par
définition de $Q_\varepsilon$, pour tout $\omega\in
\mathcal H_0$, $Q_\varepsilon(\omega)$ tend vers $Q(\omega)=
\int_U|\omega|^2\de v_g$, c'est-à-dire la norme $L^2$ au carré de $\omega$
restreinte à $U$ et $Q_\varepsilon\geq Q$ pour tout $\varepsilon$.

Passons à l'hypothèse (ii).
 Par définition de $\mathcal H_1$ et $\mathcal H_2$, il existe
pour tout $\omega\in\mathcal H_0$ un prolongement $\tilde
\varphi_\omega$ de $\varphi_\omega$ tel que $\de\tilde\varphi_\omega=\omega$.
On a alors
\begin{equation}
|\omega|_\varepsilon^2\leq\|\tilde\varphi_\omega\|_\varepsilon^2\to
\int_U|\tilde\varphi_\omega|^2\de v_g=|\omega|^2.
\end{equation}
En outre, si $\de\varphi=\omega$, alors $\varphi$ et $\tilde\varphi_\omega$ 
ne diffèrent que par une forme fermée de $M$, et donc leur restriction à 
$\overline U$ ne diffèrent aussi que par une forme fermée. Par conséquent, la
norme de $\varphi_\omega$ minore la norme sur $U$ de $\varphi$, et 
\begin{equation}
|\omega|^2=\int_U|\tilde\varphi_\omega|^2\de v_g\leq\int_U|\varphi|^2\de v_g
\leq\|\varphi\|_\varepsilon^2
\end{equation}
On en déduit que $|\omega|\leq|\omega|_\varepsilon$, et donc 
que $|\omega|_\varepsilon\to|\omega|$ pour tout $\omega$.

L'hypothèse (i), dans sa version faible (cf. remarque~\ref{conv:rem1}),
est la plus technique à vérifier.
On doit contrôler $|\cdot|_\varepsilon$ en fonction de $|\cdot|$ et
$Q$. Pour ce faire, on fixe un élément $\omega\in\mathcal H_0$ et on va 
construire une primitive particulière dont la norme $L^2$ majorera 
$|\omega|_\varepsilon$.

 On note $\bar\varphi$ la $p$-forme définie par $\bar\varphi=\varphi_\omega$
sur $\overline U$, et prolongée harmoniquement (pour la métrique $g$) sur $M$. 
On a alors $\|\bar\varphi\|_{H^1(M\backslash U)}\leq C 
\|\varphi_{|\partial U}\|_{H^\frac12(\partial U)}$, 
la constante $C$ étant indépendante de $\varphi$ (la démonstration
est similaire au cas des fonctions, cf.~\cite{ta96}, ch.~5, prop.~1.7). 
Par définition de 
$\mathcal H_0$, la $(p+1)$-forme $\bar\omega_\infty=\de\bar\varphi-\omega$ 
est un élément de $\mathcal H_\infty$ et elle est orthogonale à $\omega$. 
Elle vérifie donc $\|\bar\omega_\infty\|^2_{L^2(M\backslash U)}
\leq\|\de\bar\varphi\|^2_{L^2(M\backslash U)}\leq
\|\bar\varphi\|_{H^1(M\backslash U)}^2$ et 
elle admet une primitive $\bar\varphi_\infty$ nulle sur $U$ dont la norme
vérifie $\|\bar\varphi_\infty\|^2\leq\|\bar\omega_\infty\|^2/
\lambda^{(D)}$. Si on pose $\varphi=\bar\varphi-\bar\varphi_\infty$,
on a alors $\de\varphi=\omega$ et 
\begin{eqnarray}
\|\varphi\|_{L^2(M\backslash U)}&\leq&\|\bar\varphi\|_{L^2(M\backslash U)}+
\|\bar\varphi_\infty\|_{L^2(M\backslash U)}\nonumber\\
&\leq& \|\bar\varphi\|_{L^2(M\backslash U)}+
\|\bar\varphi\|_{H^1(M\backslash U)} /{\lambda^{(D)}}^{1/2},
\end{eqnarray}
Toutes les normes étant ici relative à la métrique $g$.

Comme $\|\bar\varphi\|_{H^1(M\backslash U)}\leq C 
\|\varphi_{|\partial U}\|_{H^\frac12(\partial U)}$ et que la norme 
$\|\varphi_{|\partial U}\|_{H^\frac12(\partial U)}$ est elle-même contrôlée 
par la norme $H^1$ de $\varphi$ sur $U$, on a finalement
\begin{equation}\label{conv:eq}
\int_{M\backslash U}|\varphi|^2\de v_g\leq C'\|\varphi\|_{H^1(U)},
\end{equation}
où $C'$ est une constante dépendant de $g$ mais pas de $\varepsilon$.

En utilisant le fait que $\varphi$ est cofermée sur $U$ (car égale à 
$\varphi_\omega$) et tangentielle le long de $\partial U$, une inégalité 
elliptique à trace associée à l'opérateur $\de+\codiff$ (voir~\cite{ta96}, 
section~5.9) donne 
\begin{equation}\label{conv:eq2}
\|\varphi\|_{H^1(U)}
\leq C''(\|\varphi\|_{L^2(U)}+\|\de\varphi\|_{L^2(U)})=
C''(|\omega|+Q(\omega)^\frac12),
\end{equation}
la métrique considérée étant ici encore $g$.

Pour une métrique $g_\varepsilon$, on a $|\omega|_\varepsilon^2\leq
\|\varphi\|_{L^2(U)}^2+\varepsilon^{n-2p}\|\varphi\|_{L^2(M\backslash U)}^2$.
Comme $\|\varphi\|_{L^2(U)}=|\omega|$ et que $\|\varphi\|_{L^2(M\backslash U)}$
peut être majoré à l'aide de~(\ref{conv:eq}) et~(\ref{conv:eq2}), 
on obtient la majoration 
\begin{equation}
|\omega|_\varepsilon^2\leq |\omega|^2+\varepsilon^{n-2p}C'C''^2
(|\omega|+Q(\omega)^\frac12)^2
\end{equation}
qui permet d'appliquer la remarque~\ref{conv:rem1} et le lemme~\ref{conv:lem2}.
\end{demo}
\section{Prescription du spectre}\label{presc}
Pour construire des valeurs propres multiples  nous allons nous utiliser,
outre théorème de convergence spectrale~\ref{conv:th}, une propriété
de transversalité vérifée par des valeurs propres multiples sur des
espaces modèles. Cette propriété remonte à Arnol'd et a été précisée
par Y.~Colin de Verdière dans \cite{cdv88}, nous allons en rappeler
une définition: 

On suppose qu'on a une famille d'opérateurs $(P_a)_{a\in B^k}$, où $B^k$
est la boule unité de $\R^k$ (en pratique, $P_a$ est le laplacien
associé à une métrique $g_a$), tels que
$P_0$ possède une valeur propre $\lambda_0$ d'espace propre $E_0$ et 
de multiplicité $N$. Pour les petites valeurs de $a$, $P_a$ possède des
valeurs propres proches de $\lambda_0$ dont la somme des espaces propres
est de dimension~$N$. Comme dans la définition de l'écart spectral,
on identifie cette somme à $E_0$ et on note $q_a$ la forme
quadratique associée à $P_a$ transportée sur $E_0$.
\begin{definition}\label{presc:def}
On dit que $\lambda_0$ vérifie l'hypothèse de transversalité d'Arnol'd si 
l'application $\Psi:a\mapsto q_a$ de $B^k$ dans $\mathcal Q(E_0)$ est 
essentielle en $0$, c'est-à-dire qu'il existe $\varepsilon>0$ tel que si 
$\Phi:B^k\to\mathcal Q(E_0)$ vérifie $\|\Psi-\Phi\|_{\infty}\leq\varepsilon$, 
alors il existe $a_0\in B^k$ tel que $\Phi(a_0)=q_0$.
\end{definition}
Une propriété cruciale est que si $\Phi$ provient d'une famille $(P'_a)$
d'opérateurs, alors $\lambda_0$ est valeur propre de $P'_{a_0}$ de 
multiplicité~$N$ et vérifie la même propriété de transversalité : on dit
que cette valeur propre multiple est stable. Comme remarqué dans \cite{cdv88},
on peut généraliser cette définition à une suite finie de valeur propre.

Pour démontrer les théorèmes~\ref{intro:th1} et~\ref{intro:th2}, nous
allons construire des domaines modèles dont le début du spectre 
vérifie la propriété de stabilité pour ensuite appliquer le 
théorème~\ref{conv:th}. Ces domaines seront des produits d'une sphère et
d'une boule, ils pourront donc être plongé dans n'importe quelle variété
de même dimension.

 La multiplicité stable sur un produit sera obtenue grâce à la formule
de Künneth : si $(M_1,g_1)$ et $(M_2,g_2)$ sont deux variétés
riemanniennes compactes et $\alpha_i\in\Omega^*(M_i)$, $i=1,2$, alors
on a, en identifiant chacune des formes $\alpha_i$ à son relevé sur
$(M_1\times M_2,g_1\oplus g_2)$,
\begin{equation}
\Delta(\alpha_1\wedge\alpha_2)=\Delta\alpha_1\wedge\alpha_2+\alpha_1
\wedge\Delta\alpha_2.
\end{equation}
En particulier, si $\alpha_i$ est une forme propre de valeur propre $\lambda_i$
pour $i=1,2$, alors $\alpha_1\wedge\alpha_2$ est une forme propre de valeur
propre $\lambda_1+\lambda_2$ pour la métrique produit. Il est clair que
que si $\alpha_1$ et $\alpha_2$ sont fermées, alors $\alpha_1\wedge\alpha_2$
aussi. On peut vérifier que pour la métrique produit, si $\alpha_1$ et 
$\alpha_2$ sont cofermées (par exemple si l'une des deux est une fonction)
leur produit est aussi cofermé. 

On va montrer que si l'une des deux valeurs propres est de multiplicité 
stable, la valeur propre multiple obtenue sur la variété produit peut 
hériter de cette propriété.
\begin{lem}
On suppose que $\mu_{p_1,k_1}(M_1,g_1)$ est une valeur propre
simple, $\mu_{p_2,k_2}(M_2,g_2)$ une valeur propre de multiplicité stable $N$
et qu'il n'y a pas d'autres valeurs propres de $M_1$ et $M_2$ dont 
la somme soit $\mu_{p_1,k_1}(M_1,g_1)+\mu_{p_2,k_2}(M_2,g_2)$.
Alors $\mu_{p_1,k_1}(M_1,g_1)+\mu_{p_2,k_2}(M_2,g_2)$ est valeur
propre de multiplicité stable $N$ sur $(M_1\times M_2,g_1\oplus g_2)$.
\end{lem}
\begin{demo}
Si on note $\alpha$ une forme propre de $\mu_{p_1,k_1}(M_1,g_1)$ et 
$E$ l'espace propre de $\mu_{p_2,k_2}(M_2,g_2)$ alors, en vertu des 
remarques qui précèdent, $\alpha\wedge E$ est un espace propre cofermé de 
valeur propre $\mu_{p_1,k_1}(M_1,g_1)+\mu_{p_2,k_2}(M_2,g_2)$.

Supposons que $\mu_{p_2,k_2}(M_2,g_2)$ vérifie la propriété~\ref{presc:def}
pour la famille de métrique $g_{2,a}$ et la famille de forme quadratique
$q_a$. Si la forme $\alpha$ est normée, le produit par 
$\alpha$ plonge isométriquement $\Omega(M_2)$ dans $\Omega(M_1\times M_2)$
en envoyant espace propre sur espace propre, quel que soit $a$.
En particulier, les espaces propres de $(M_2,g_{2,a})$ de valeur propre 
proche de $\mu_{p_2,k_2}(M_2,g_2)$, sont envoyé sur des espaces propres
de même multiplicité, leur valeur propre étant augmenté de 
$\mu_{p_1,k_1}(M_1,g_1)$. L'espace $\alpha\wedge E$ vérifie donc de la 
propriété~\ref{presc:def} pour la valeur 
propre $\mu_{p_1,k_1}(M_1,g_1)+\mu_{p_2,k_2}(M_2,g_2)$  et la famille de 
forme quadratique $q_a+\mu_{p_1,k_1}(M_1,g_1)|\cdot|^2$.
\end{demo}

Nous utiliserons plusieurs fois le fait que sur une boule euclidienne
de rayon $\varepsilon$, le noyau du laplacien pour la condition de bord (A)
est trivial, sauf en degré 0 pour lequel sa dimension est 1. Toutes les
autres valeurs propres tendent vers $+\infty$ quand $\varepsilon\to0$. En 
particulier, sur un produit riemannien $M\times B^k(\varepsilon)$, les
premières formes propres sont les relevés des premières formes propres de $M$.

Les trois lemmes qui suivent ont pour but de construire les domaines modèles:
\begin{lem}\label{presc:lem1}
Pour tous entiers $N\geq1$, $n\geq3$, toute suite finie $0<a_1\leq a_2\leq
\ldots\leq a_N$ et toute constante $C>a_N$, il existe une métrique $g$ sur 
$S^n$ telle que $\mu_{0,i}=a_i$ pour $i\leq N$, ces valeurs propres 
vérifiant l'hypothèse de stabilité, $\mu_{0,N+1}>C$ et
$\mu_{p,1}>C$ pour $1\leq p\leq[\frac{n-1}2]$. 
\end{lem}
\begin{demo}
Le résultat sur les $\mu_{0,i}$ découle des travaux de Colin de Verdière
(voir \cite{cdv86} et \cite{cdv87}), il suffit de montrer que la construction
géométrique peut se faire avec $\mu_{p,1}>C$ pour $p\geq1$. Le principe de 
la démonstration
est d'appliquer le théorème~\ref{conv:th} avec un domaine $U$ dont le spectre
vérifie la conclusion du lemme. Pour appliquer l'argument de stabilité,
on ne doit pas travailler avec une seule métrique sur $U$ mais une petite
famille de métriques, un élément crucial sera alors que l'invariant 
conforme qui minore le spectre en degré $[\frac n2]$ varie continûment
pour la topologie $C^0$ (voir remarque~\ref{conv:rq2}), il sera donc 
uniformément minoré pour cette famille de métrique.

On procède par récurrence sur la dimension. En dimension 3 et 4, on procède
comme dans \cite{cdv87}: On choisit une surface $\Sigma$ dont le début 
du spectre (pour les fonctions) est égal à $a_1\leq a_2\leq\ldots\leq a_N$ 
et vérifie l'hypothèse de stabilité, et on choisit comme domaine $U$ 
le produit riemannien de $\Sigma$ avec un petit intervalle 
$]-\varepsilon,\varepsilon[$ (en dimension~3) ou un petit disque de rayon 
$\varepsilon$ (en dimension~4), la métrique étant prolongée de manière 
quelconque en dehors de $U$. Le théorème~\ref{conv:th} et l'argument de 
stabilité assure l'existence d'une métrique $g$ telle que $\mu_{0,i}=a_i$ 
pour $i\leq N$, et $\mu_{0,N+1}>C$, et comme son volume est arbitrairement 
petit, le théorème~\ref{conv:th2} et la remarque~\ref{conv:rq2} assurent 
qu'on peut choisir $g$ telle qu'on ait aussi $\mu_{1,1}>C$.

 Si, par hypothèse de récurrence, le lemme est vrai en dimension $n-1$, 
on raisonne comme en dimension~3 et~4 en munissant $S^{n-1}$ de la métrique
fournie par le lemme et en utilisant le domaine $U=S^{n-1}\times
]-\varepsilon,\varepsilon[$ plongé dans $S^n$. Le théorème~\ref{conv:th}
donne une métrique $S^n$ qui a le spectre souhaité pour $p\leq[\frac{n-3}2]$
($H^p(U/M)$ étant trivial), et le théorème~\ref{conv:th2} assure que 
$\mu_{[\frac{n-1}2],1}>C$.
\end{demo}

\begin{lem}\label{presc:lem2}
Pour tous entiers $N\geq1$, $p\geq2$ et $n\geq 2p+3$, toute suite finie 
$0<a_1\leq a_2\leq\ldots\leq a_N$ et toute constante $C>a_N$, il existe 
une métrique $g$ sur la variété $M^n=S^{p+1}\times B^{p+2}$ pour $n=2p+3$
ou $M^n=S^{n-1}\times [0,1]$ si $n>2p+3$, telle que $\mu_{p,i}=a_i$ pour 
$i\leq N$, ces valeurs propres vérifiant l'hypothèse de stabilité, 
$\mu_{p,N+1}>C$ et $\mu_{q,1}>C$ pour $1\leq q\leq\frac n2$, 
$q\neq p$, le volume $\Vol(M,g)$ étant arbitrairement petit.
\end{lem}
\begin{demo}
Remarquons d'abord qu'en prescrivant le début du spectre pour $p=0$
dans le lemme~\ref{presc:lem1}, on a aussi prescrit les $\mu_{p,i}$
pour $p=n-1$ et $i\leq N$, par dualité de Hodge. Partant de ce constat,
on va encore procéder par récurrence sur $n$, le degré $p$ étant fixé.

Pour $n=2p+3$, il suffit de considérer $M=S^{p+1}\times B^{p+2}(\varepsilon)$,
la sphère $S^{p+1}$ étant munie de la métrique fournie par le 
lemme~\ref{presc:lem1} et $B^{p+2}(\varepsilon)$ étant une boule de rayon 
$\varepsilon$ petit.

Pour $n>2p+3$, la récurrence s'effectue comme dans le lemme~\ref{presc:lem1}:
on applique le théorème~\ref{conv:th} au domaine $M^{n-1}$ plongé dans 
$S^{n-1}$ et on pose $M^n=S^{n-1}\times]-\varepsilon,\varepsilon[$.
\end{demo}

\begin{lem}\label{presc:lem3}
Pour tous entiers $N\geq1$, $n\geq 3$, toute suite finie
$0<a_1<a_2\leq a_3\leq\ldots\leq a_N$ et toute constante $C>a_N$, il existe
une métrique $g$ sur $M=B^3\times S^n$ telle que $\mu_{1,i}(M,g)=a_i$ pour 
$i\leq N$ (avec stabilité), $\mu_{1,N+1}(M,g)>C$ et $\mu_{p,1}(M,g)>C$ pour 
$2\leq p\leq\frac n2$, le volume $\Vol(M,g)$ étant arbitrairement petit.
\end{lem}
\begin{demo}
On va une nouvelle fois utiliser le théorème~\ref{conv:th} avec un domaine
$U$ produit d'une sphère et d'une boule, mais avec une métrique particulière
sur la boule : selon \cite{gu04} (théorème~2.1), il existe une métrique $g_B$ 
sur $B^3$ telle que $\mu_{1,1}(B^3,g_B)=a_1$, $\mu_{1,2}(B^3,g_B)>C$ et 
$\mu_{p,1}(B^3,g_B)>C$ pour $p=0,2$, le volume étant arbitrairement petit. 
En utilisant le lemme~\ref{presc:lem1}, on peut munir $S^n$ d'une métrique 
telle que $\mu_{0,i}(S^n)=a_i-a_1$ pour $i\leq N-1$, $\mu_{0,N}(S^n)>C$
et $\mu_{p,1}(S^n)>C$ pour $p\geq1$. Le produit $M=S^n\times B^3$
vérifie alors $\mu_{1,i}(M)=\mu_{0,i}(S^n)+\mu_{1,1}(B^3,g_B)=a_i$
pour $i\leq N$, $\mu_{1,N+1}(M)>C$ et $\mu_{p,1}(M)>C$ pour 
$2\leq p\leq\frac n2$. 
\end{demo}

On est maintenant en mesure de démontrer les résultats annoncés dans
l'introduction.

\begin{demo}[du théorème~\ref{intro:th1}]
Pour chaque degré $p$, on se donne le domaine $U_p$ fourni par le 
lemme~\ref{presc:lem2} (ou le lemme~\ref{presc:lem3} pour $p=1$) dont
le début du spectre en degré $p$ est celui qu'on veut prescrire, avec 
$C>\sup_{p,i}a_{p,i}$. On peut noter que dans la démonstration du 
théorème~\ref{conv:th}, l'hypothèse de connexité du domaine $U$ 
n'intervient pas, la conclusion reste donc vraie en remplaçant $U$ par 
un nombre fini de domaine. On peut donc l'appliquer à la famille $(U_p)$
plongée dans $M$ (en remarquant que $H^p((\cup U_p)/M)$ est trivial pour
$1\leq p\leq\frac{n-3}2$), et grâce à la stabilité du spectre de ces domaines
on conclut à l'existence d'une métrique $g$ sur $M$ telle que
$\mu_{p,k}(M,g)=a_{p,k}$ pour pour $1\leq k\leq N$ et 
$1\leq p\leq[\frac{n-3}2]$. Comme pour les lemmes précédents, le cas du 
degré $p=[\frac{n-1}2]$ est couvert par le théorème~\ref{conv:th2}. 

Pour montrer qu'on peut aussi prescrire le volume, on procède comme dans
\cite{gu04} et \cite{ja08} : on peut appliquer l'argument de stabilité
au volume en le traitant comme une valeur propre simple. Il suffit donc
d'ajouter à la famille $U_p$ une boule de volume $V-\sum_p\Vol(U_p)$ 
(en ayant choisit les volumes des $U_p$ suffisamment petits) et dont 
les valeurs propres pour les $p$-formes sont arbitrairement grandes
(c'est possible selon \cite{gp95}). Le fait que le volume de cette boule
vérifie l'hypothèse de transversalité signifie simplement qu'on
peut lui donner n'importe quelle valeur au voisinage de $V-\sum_p\Vol(U_p)$,
par exemple par homothétie.
\end{demo}

\begin{demo}[du théorème~\ref{intro:th2}]
Y.~Colin de Verdière a montré que la première valeur propre de la sphère
$S^2$ muni de sa métrique canonique est stable (\cite{cdv88}, section~2).
On peut donc raisonner comme dans le lemme~\ref{presc:lem2}: la
valeur propre $\mu_{1,1}(S^2)$ est de multiplicité~3 stable et il suffit
d'appliquer le théorème~\ref{conv:th} au domaine 
$S^2\times B^{n-2}(\varepsilon)$.
\end{demo}

\noindent Pierre \textsc{Jammes}\\
Université d'Avignon et des pays de Vaucluse\\
Laboratoire d'analyse non linéaire et géométrie (EA 2151)\\
F-84018 Avignon\\
\texttt{Pierre.Jammes@ens-lyon.org}

\begin{thebibliography}{HHN99}
{\small
\makeatletter
\ifx\fonteauteurs\@undefined
\def\fonteauteurs{\scshape}\fi
\makeatother

\bibitem[Am03]{am03}
\bgroup\fonteauteurs B.~Ammann\egroup{} -- \og~A spin-conformal lower bound of
  the first positive {D}irac eigenvalue~\fg, {\em Differ. Geom. Appl.}, 18 (1),
  p.~21--32, 2003.

\bibitem[An89]{an89}
\bgroup\fonteauteurs C.~Anné\egroup{} -- \og~Principe de {D}irichlet pour les
  formes différentielles~\fg, {\em Bull. soc. math. France}, 117 (4),
  p.~445--450, 1989.

\bibitem[BCC98]{bcc98}
\bgroup\fonteauteurs G.~Besson\egroup{}, \bgroup\fonteauteurs
  B.~Colbois\egroup{} et \bgroup\fonteauteurs G.~Courtois\egroup{} -- \og~Sur
  la multiplicité de la première valeur propre de l'opérateur de {S}chrödinger
  avec champ magnétique sur la sphère ${S}^2$~\fg, {\em Trans. Amer. Math.
  Soc.}, 350 (1), p.~331--345, 1998.

\bibitem[BD97]{bd97}
\bgroup\fonteauteurs G.~Baker\egroup{} et \bgroup\fonteauteurs
  J.~Dodziuk\egroup{} -- \og~Stability of spectra of {H}odge-de~{R}ham
  laplacians~\fg, {\em Math. Z.}, 224 (3), p.~327--345, 1997.

\bibitem[Be80]{be80}
\bgroup\fonteauteurs G.~Besson\egroup{} -- \og~Sur la multiplicité de la
  première valeur propre des surfaces riemanniennes~\fg, {\em Ann. inst.
  Fourier}, 30 (1), p.~109--128, 1980.

\bibitem[CdV86]{cdv86}
\bgroup\fonteauteurs Y.~Colin~de Verdière\egroup{} -- \og~Sur la multiplicité
  de la première valeur propre non nulle du laplacien~\fg, {\em Comment. Math.
  Helv.}, 61 (2), p.~254--270, 1986.

\bibitem[CdV87]{cdv87}
\bgroup\fonteauteurs Y.~Colin~de Verdière\egroup{} -- \og~Construction de
  laplaciens dont une partie finie du spectre est donnée~\fg, {\em Ann. scient.
  \'Ec. norm. sup.}, 20 (4), p.~99--615, 1987.

\bibitem[CdV88]{cdv88}
\bgroup\fonteauteurs Y.~Colin~de Verdière\egroup{} -- \og~Sur une hypothèse de
  transversalité d'{A}rnol'd~\fg, {\em Comment. Math. Helv.}, 63 (2),
  p.~184--193, 1988.

\bibitem[CdVT93]{cdvt93}
\bgroup\fonteauteurs Y.~Colin~de Verdière\egroup{} et \bgroup\fonteauteurs
  N.~Torki\egroup{} -- \og~Opérateur de {S}chrödinger avec champ
  magnétique~\fg, {\em Sémin. théor. spectr. géom.}, 11, p.~9--18, 1993.

\bibitem[Ch76]{ch76}
\bgroup\fonteauteurs S.~Y. Cheng\egroup{} -- \og~Eigenfunctions and nodal
  sets~\fg, {\em Comment. Math. Helv.}, 51 (1), p.~43--55, 1976.

\bibitem[Ch79]{ch79}
\bgroup\fonteauteurs J.~Cheeger\egroup{} -- \og~Analytic torsion and the heat
  equation~\fg, {\em Ann. Math.}, 109 (2), p.~259--322, 1979.

\bibitem[Co04]{co04}
\bgroup\fonteauteurs B.~Colbois\egroup{} -- \og~Spectre conforme et métriques
  extrémales~\fg, {\em Sémin. théor. spectr. géom.}, 22, p.~93--101, 2004.

\bibitem[Da05]{da05}
\bgroup\fonteauteurs M.~Dahl\egroup{} -- \og~Prescribing eigenvalues of the
  {D}irac operator~\fg, {\em Manuscripta Math.}, 118 (2), p.~191--199, 2005,
  math.DG/0311172.

\bibitem[Do82]{do82}
\bgroup\fonteauteurs J.~Dodziuk\egroup{} -- \og~Eigenvalues of the {L}aplacian
  on forms~\fg, {\em Proc. of Am. Math. Soc.}, 85, p.~438--443, 1982.

\bibitem[Er02]{er02}
\bgroup\fonteauteurs L~Erd{\H{o}}s\egroup{} -- \og~Spectral shift and
  multiplicity of the first eigenvalue of the magnetic {S}chr\"odinger operator
  in two dimensions~\fg, {\em Ann. inst. Fourier}, 52 (6), p.~1833--1874, 2002.

\bibitem[GP95]{gp95}
\bgroup\fonteauteurs G.~Gentile\egroup{} et \bgroup\fonteauteurs
  V.~Pagliara\egroup{} -- \og~Riemannian metrics with large first eigenvalue on
  forms of degree $p$~\fg, {\em Proc. of Am. Math. Soc.}, 123 (12),
  p.~3855--3858, 1995.

\bibitem[Gu04]{gu04}
\bgroup\fonteauteurs P.~Guérini\egroup{} -- \og~Prescription du spectre du
  laplacien de {H}odge-de~{R}ham~\fg, {\em Ann. scient. \'Ec. norm. sup.}, 37
  (2), p.~270--303, 2004.

\bibitem[HHN99]{hohon99}
\bgroup\fonteauteurs M.~Hoffmann-Ostenhof\egroup{}, \bgroup\fonteauteurs
  T.~Hoffmann-Ostenhof\egroup{} et \bgroup\fonteauteurs
  N.~Nadirashvili\egroup{} -- \og~On the multiplicity of eigenvalues of the
  {L}aplacian on surfaces~\fg, {\em Ann. Global Anal. Geom.}, 17 (1),
  p.~43--48, 1999.

\bibitem[Ja07]{ja07}
\bgroup\fonteauteurs P.~Jammes\egroup{} -- \og~Minoration conforme du spectre
  du laplacien de {H}odge-de~{R}ham~\fg, {\em Manuscripta Math.}, 123 (1),
  p.~15--23, 2007, math.DG/0604591.

\bibitem[Ja08]{ja08}
\bgroup\fonteauteurs P.~Jammes\egroup{} -- \og~Prescription du spectre du
  laplacien de {H}odge-de~{R}ham dans une classe conforme~\fg, {\em Comment.
  Math. Helv.}, 83 (3), p.~521--537, 2008, math.DG/0601738.

\bibitem[Ja09a]{ja09a}
\bgroup\fonteauteurs P.~Jammes\egroup{} -- \og~Construction de valeurs propres
  doubles du laplacien de {H}odge-de~{R}ham~\fg, {\em J. Geom. Anal.}, 19 (3),
  p.~643--654, 2009, math.DG/0608758.

\bibitem[Ja09b]{ja09b}
\bgroup\fonteauteurs P.~Jammes\egroup{} -- \og~Sur la multiplicité des valeurs
  propres d'une variété compacte~\fg, {\em S\'emin. th\'eor. spectr. g\'eom.},
  26, p.~1--11, 2009.

\bibitem[Lo96]{lo96}
\bgroup\fonteauteurs J.~Lohkamp\egroup{} -- \og~Discontinuity of geometric
  expansions~\fg, {\em Comment. Math. Helv.}, 71 (2), p.~213--228, 1996.

\bibitem[{Mc}93]{mc93}
\bgroup\fonteauteurs J.~{Mc}{Gowan}\egroup{} -- \og~The $p$-spectrum of the
  {L}aplacian on compact hyperbolic three manifolds~\fg, {\em Math. Ann.}, 297
  (4), p.~725--745, 1993.

\bibitem[Na88]{na88}
\bgroup\fonteauteurs N.~Nadirashvili\egroup{} -- \og~Multiple eigenvalues of
  the {L}aplace operator~\fg, {\em Math. USSR-Sb.}, 61 (1), p.~225--238, 1988.

\bibitem[Sé94]{se94}
\bgroup\fonteauteurs B.~Sévennec\egroup{} -- \og~Multiplicité du spectre des
  surfaces: une approche topologique~\fg, {\em Sémin. théor. spectr. géom.},
  12, p.~29--36, 1994.

\bibitem[Sé02]{se02}
\bgroup\fonteauteurs B.~Sévennec\egroup{} -- \og~Multiplicity of the second
  {S}chrödinger eigenvalue on closed surfaces~\fg, {\em Math. Ann.}, 324 (1),
  p.~195--211, 2002.

\bibitem[Ta96]{ta96}
\bgroup\fonteauteurs M.~Taylor\egroup{} -- {\em Partial differential equations
  I}, Springer, 1996.

}\end{thebibliography}
\end{document}